\documentstyle[12pt,righttag]{amsart}
\newtheorem{prop}{Proposition}
\newtheorem{th}[prop]{Theorem}

\newtheorem{lem}[prop]{Lemma}

\theoremstyle{definition}
\newtheorem{ack}{Acknowledgments} 
\theoremstyle{remark}
\newtheorem{rem}{Remark}

\newcommand{\bbR}{{\Bbb{R}}}

\newcommand{\s}{\sigma}

\newcommand{\Om}{\Omega}
\newcommand{\La}{\Lambda}

\newcommand{\G}{\Gamma}

\renewcommand{\span}{\operatorname{span}}
\newcommand{\supp}{\operatorname{supp}}

\newcommand{\dint}{\displaystyle\int}

\newcommand{\disp}{\displaystyle}
\newcommand{\lb}{\label}

\newcommand{\lra}{\longrightarrow}

\makeatletter
\def\@currentlabel{2.1}\label{e:dispaa}
\def\@currentlabel{2.21}\label{e:dispau}
\def\@currentlabel{2.22}\label{e:dispav}
 \def\@currentlabel{2.23}\label{e:dispaw}
\def\@currentlabel{2.24}\label{e:dispax}
\makeatother

\makeatletter
\def\alphenumi{%
  \def\theenumi{\alph{enumi}}%
  \def\p@enumi{\theenumi}%
  \def\labelenumi{(\@alph\c@enumi)}}
\makeatother

\begin{document}

\title{Contractive projections in nonatomic function spaces}
\author{Beata Randrianantoanina}
\address{Department of Mathematics \\ Bowling Green State University \\
         Bowling Green, OH 43403}
\email{brandri@@andy.bgsu.edu}
\maketitle

\begin{abstract}
We prove that there is no 1-complemented subspaces of finite
codimension in separable rearrangament-invariant function spaces.
\end{abstract}

We study contractive  projections onto
finite-codimensional subspaces of real nonatomic function spaces. In
general such
projections are not common. It is well known that only in Hilbert space
there exists a contractive projection onto every subspace of fixed finite
codimension (cf.~\cite{Amir}).

  The study of contractive projections (and more general projections with
minimal norm) is important in approximation theory (cf. the
survey of Cheney and Price~\cite{ChP}).

  It is known  that there is no 1-com\-ple\-ment\-ed
subspaces of finite-codimension in
$C[0,1]$ (Wulbert \cite{Wbt}) and
    in $L_1(\mu)$ if the measure $\mu$ is
nonatomic (\cite[Corollary IV.1.15]{HWW},  \cite{D}).
De Figueiredo and Karlovitz \cite{FK} (1970) proved that if $\mu$ is
nonatomic then
there is no 1-complemented hyperplanes
in $L_p(\Om,\mu),\  1<p<\infty,\ p\neq 2,$
 (cf. also \cite{BM}).

In this paper we prove
that in rearrangement-invariant nonatomic function
spaces not isometric to $L_2$ there is no 1-complemented subspaces of any
finite codimension.

We use the terminology and notation as in \cite{LT2}.

Our method of proof is surprisingly simple -- it is based on the following
observation:

\begin{prop}  \label{A} $(\text{cf.} $ \cite{KR,R84}$)$
 In a real Banach space $X$ if    $P$ is a projection
 then $\|I - P\|=1 $ (where $I$ denotes identity operator)
 if and only if $x^*(Px) \ge 0$ for all $x\in X$ and $x^* \in X^*$
norming for $x.$
\end{prop}

In \cite[Theorem 4.3]{KR} (cf. \cite{KR93,th})   Kalton and the
author proved the nonexistence of 1-complemented hyperplanes in a wide
class of
nonatomic function spaces. However the original theorem uses
special technical frazeology so we state it below in the language
of projections:
\begin{th}  \lb{kw}
Suppose
$X$ is a real order-continuous  K\"othe function space on
$(\Omega,\mu)$ and $\mu$ is nonatomic.
Then the hyperplane $H$ in $X$ is 1-complemented
 if and only if there exists a nonnegative measurable
function $w $ with $\supp w=B= \supp f,$ where $f\in H^\perp
\subset X^*,$ so that for any $x\in X$ with $\supp x
\subset B$
 $$  \|x\|
=(\int
|x|^2w\,d\mu)^{1/2}.$$
\end{th}

Hence there is no 1-complemented hyperplanes in $X$ unless
$L_2$ is isometric to a band in $X.$ In particular there is no
1-complemented hyperplanes in separable r.i. spaces on
$[0,1]$ (\cite[Theorem 4.4]{KR}).

As a part of proof of Theorem~\ref{kw}
 we proved the following fact which we state
separately for the future use.

\begin{prop} $($\cite[Proposition 2.8]{th}$)$ \lb{l2}
Let $X$ be an order-continuous  K\"othe function space
on $(\Omega,\mu),$ where $\mu$ is nonatomic. Suppose
that  the set $${\disp \La = \left\{\frac{x^*}{x} \ : x\in X,
x^* \in X^* \ \text{norming \ for } x \right\} }$$ is
one-dimensional i.e. $\La \subset \{ aw  : a \in \bbR\}$
for some  $w \in L_0(\Om,\mu).$   Then $X$ is isometric
to $L_2(wd\mu).$
\end{prop}

Now we are ready to prove our main result.

\begin{th} \lb{fc}
  Suppose $\mu$ is nonatomic and   $X$ is a separable r.i.
 space on
$([0,1],\mu)$ not isometric to $L_2.$ Then there is no
1-complemented   subspaces of any
finite codimension in $X.$
\end{th}

For the proof   we   need the following
measure-theoretic lemma.

\begin{lem}   \lb{liap2}
  Suppose $\mu$ is nonatomic and suppose $f_1,\ldots ,f_n,g_1,\ldots ,g_n
\in L_1(\mu)$ are such that $g_1,\ldots ,g_n$ are linearly independent
and
 $$  \sum_{i=1}^n \left( \int hf_j\ d\mu \right)\left( \int hg_j\ d\mu
                                                 \right) \ge 0
 $$
whenever $|h| = 1$~a.e. Then  $\{f_j\}_{j=1}^n \subset
\span{\{g_j\}_{j=1}^n  }.$
\end{lem}

\begin{pf}
Consider the subset $\G$ of
$\bbR^{2n}$ of all $2n$-tuples $(a_1,b_1,\ldots , a_n,b_n)$ such that for
some $h \in L_{\infty}(\mu)$ with $|h| = 1$~a.e. we have
$a_j = \dint hf_j\ d\mu,\ \ b_j = \dint hg_j\ d\mu,\ j=1,\ldots ,n.$
Let $H = \span{\G} \subset \bbR^{2n}.$

We immediately see that $\G = -\G$ and that by Liapunoff's theorem
(\cite{Rud}) $\G$ is convex. Hence for every $(s_1,t_1,\ldots , s_n,t_n)
\in H$ we have ${\disp\sum_{j=1}^n s_j t_j \ge 0.}$

On the other hand consider
$V = \left\{ v_k = e_{2k-1} - e_{2k}\/:k= 1,\ldots ,n \right\}
      \subset \bbR^{2n}$
where $e_j$ denotes the natural basis of $\bbR^{2n}.$ Clearly
$V\cap H = \{ 0\} $ so $\dim H \le n.$
Therefore there exist $c_1, \ldots ,c_n \in \bbR$ such that
$\disp{ s_i = \sum_{j=1}^n c_jt_j}$ for $i = 1, \ldots ,n,$\/
for every
$(s_1,t_1,\ldots , s_n,t_n) \in H.$
In particular
${\disp \int hf_i\ d\mu \/=\/ \int h\left( \sum_{j=1}^n c_jg_j \right)\/
   d\mu}$
for all $h$ with $|h|= 1$~a.e. and the lemma follows.
\end{pf}

\begin{pf*}{Proof of Theorem~\ref{fc}}
    Suppose that $F$ is a closed linear subspace of codimension $n$ in
$X$ and let
$u_1, \dots ,u_n$ be linearly independent functions in $X$ such that
$\span \{ F,u_1, \dots ,u_n\} = X.$
Denote $P:X\lra F$ a contractive projection onto $F$ and consider
$Q=I-P.$ Then
${\disp Q = \sum_{j=1}^n f_j\otimes u_j}$
for some linearly independent $f_1, \dots ,f_n \in X^*.$

By Proposition~\ref{A} for any $x\in X$ and $x^* \in X^*$
norming for $x$ \ $x^*(Qx)\ge 0.$
Next for any $h$ with
$|h|=1$~a.e. $hx^*$ is norming for $hx$ if $x^*$ is
norming for $x.$ Hence
 $hx^*\left( Q(hx)\right) \ge 0$ i.e.
$$ \sum_{j=1}^n \left( \int f_jhx\,d\mu \right)
                \left( \int u_jhx^*\,d\mu \right) \ge 0.$$

 By Lemma~\ref{liap2} $f_1x \in \span{\{ u_jx^*\, :j=1,\dots ,n\} }.$
So if $B = \supp{f_1} $ ($\mu(B) > 0$) then
${\disp\frac{x}{x^*}\Big|_B 
\in \span \left\{ \frac{u_j}{f_1} \, :j=1,\dots,n\right\} }.$ 
By re-arrangement invariance of $X$ for every measure preserving map
$\s : [0,1]\lra [0,1]$  $\ \ x^* \circ \s$ is norming for $x\circ \s$
and so
$${\disp\frac{x\circ \s}{x^*
\circ \s}\Big|_B  \in \span \left\{ \frac{u_j}{f_1} \, :j=1,\dots
,n\right\} }$$ i.e. the set
 ${\disp \left\{ \left[\left( \frac{x}{x^*}
 \right) \circ \s\right]\Big|_B  \, | \;
                          \s:[0,1]\lra [0,1] \text{ measure }
 \text{preserving}
\right\} } $
is finite dimensional which is impossible unless ${\disp \frac{x}{x^*} }
$ is a constant. But then by Proposition~\ref{l2} we conclude that $X$ is
isometric to $L_2[0,1]$ contrary to our assumption.
\end{pf*}

\begin{rem}
Notice that in the proof of Theorem~\ref{fc} we use
 re-arrangement invariance of $X$ only in the final step to
conclude that if the set ${\disp
\left\{\frac{x^*}{x} \ : x\in X,
x^* \in X^* \text{ norming for } x \right\} }$ is
finite-dimensional then it is one-dimensional.
Similar conclusion is true also in spaces of the form
$X(Y),$ $X_1(X_2(\ldots(X_m)\ldots)),$ where $X,Y,X_j$ are
r.i. and Theorem~\ref{fc} holds as stated also for those
spaces.
\end{rem}

\begin{ack}
I wish to express my gratitude to Professor Nigel Kalton for his interest
in this work and many valuable disscussions.
\end{ack}
 \bibliographystyle{standard}
    \bibliography{tref}

\end{document}